\documentclass[11pt,reqno]{amsart}
\usepackage{amssymb,latexsym}
\usepackage{graphicx}
\usepackage{fancyhdr,amssymb}
\usepackage{url}		%does nice formatting of URLs

\theoremstyle{plain}
\numberwithin{equation}{section}

\begin{document}
\fancyhead{}
\renewcommand{\headrulewidth}{0pt}
\fancyfoot{}
\fancyfoot[LE,RO]{\medskip \thepage}
%\fancyfoot[LO]{\medskip MONTH YEAR}
%\fancyfoot[RE]{\medskip VOLUME , NUMBER }

\setcounter{page}{1}

\title[Three term Recurrence  and residue completeness ]{Three term Recurrence
 and residue completeness }
\author{Cheng Lien Lang}
\address{Department Applied of Mathematics\\
                I-Shou University\\
                Kaohsiung, Taiwan\\
                Republic of China}
\email{cllang@isu.edu.tw}
\thanks{}
\author{Mong Lung Lang}
\address{Singapore
669608, Republic of Singapore   }
\email{lang2to46@gmail.com}

\begin{abstract}
Let $w_n= qw_{n-1}+w_{n-2}$ and let $K$ be the set of all $m$'s such that
 $\{w_n\}$ modulo $m$ is $\Bbb Z_m$. The present article studies
  $\{w_n\}$ modulo $m$  and determines the set $K$. In particular,
 (i) Pell numbers modulo $m$ is residue complete if and only if
  $m \in\{ 2, 3^a, 5^b\}$,  and (ii) Pell-Lucas numbers modulo $m$ is residue complete
 if and only if $m = 3^n$
\end{abstract}

%\medskip

\maketitle
\vspace{-.8cm}

\section {Introduction}
 Given $a, b \in \Bbb Z$, $q\in \Bbb Z  -\{0\}$.
Define the three term linear recurrence $\{w_n(a,b,q)\}=\{w_n\} $  by the following.
$$ w_0=a, w_1=b, w_{n}= qw_{n-1}+w_{n-2}.\eqno(1.1)$$
Let $m\in \Bbb N$.
 The recurrence  $\{w_n\}$ modulo $m$ is periodic (see Lemma 2.1).
 Let $w(a,b,q,m)$ be a period. The length of $w(a,b,q,m)$ is denoted by
  $k(a,b,q,m)$.
We say $\{w_n\}$ modulo $m$ is residue complete
 if
 $x\in \{w_r\}$ modulo $m$ for every $x\in \Bbb Z_m$. Equivalently,
 $x\in w(a,b,q,m)$ for every $x\in \Bbb Z_m$.
Let $p$ be a prime divisor of  such $m$.
 Applying results of Bumby [ 2], Li [6], Schinzel [7], Somer [8], and Webb and Long [9]
 $$p\in \Delta  \cup \Omega ,\eqno(1.2)$$
 where
 $\Delta = \{2,3,5,7\}$ and
 $\Omega $ is the set of prime divisors of $q^2+4$.
 Note that $3$ and $7$ cannot be divisors of $q^2+4$.
In the present article, we study $\{w_n\}$ modulo $p^e$,
where $p\in \Delta \cup \Omega$. We will give a series of lemmas that enables us to find
 the  set $K$
 that consists of all the $m$'s such that $\{w_n(a,b,q)\}$ modulo $m$ is residue complete (see Discussion 5.3).
  In particular, one has
  (i) Pell numbers modulo $m$ is residue complete if and only if
  $m \in\{ 2, 3^a, 5^b\}$,  and (ii) Pell-Lucas numbers modulo $m$ is residue complete
 if and only if $m = 3^n$ (see Proposition 5.1 and 5.2).
{\bf The period $w(0,0,q, m)$ is called the trivial period. Throughout the article, the
 periods being considered are nontrivial periods.}

 \section {Basic Properties about $w_n$}

  $$
\mbox{Let } \sigma = \left[\begin{array}{lr}
q &1 \\
1 &0\\
\end{array}\right ].
 \mbox{  Then }
  \left[\begin{array}{cr}
qb+a &b \\
b & a\\
\end{array}\right ]
\sigma^n
=
\left[\begin{array}{lr}
w_{n+2} &w_{n+1} \\
w_{n+1}&w_n\\
\end{array}\right ].
 \eqno(2.1)$$

\noindent {\bf Lemma 2.1.} {\em  The sequence $\{w_n(a,b, q)\}\,(mod\,m)$ is periodic.
Suppose  gcd$\,(a^2+qab -b^2, m)=1$. Then
 the length of
 $w(a,b,q, m)$ is
 the smallest positive integer $k$ such that
  $\sigma^k\equiv 1\,\,(mod\, m)$.
  }

\noindent {\em Proof.}
 In (2.1), there exists some $r \in \Bbb N$ such that $\sigma^r \equiv 1$ (mod $m$).
  It follows from (2.1) that $(w_1, w_0) = (b,a) \equiv (w_{r+1}, w_r)$. Hence $\{w_n\}$ modulo $m$
   is periodic.
 Let $k$ be the length of $w(a,b,q,m)$. Then $k$ is the smallest positive integer
   such that $(qb+a, b) = (w_2, w_1) \equiv  (w_{k+2}, w_{k+1})$ as well as
 $(b, a) = (w_1, w_0) \equiv (w_{k+1}, w_{k})$.
 As a consequence, identity (2.1) becomes
$$  \left[\begin{array}{cr}
qb+a &b \\
b & a\\
\end{array}\right ]
\sigma^k
=
\left[\begin{array}{lr}
w_{k+2} &w_{k+1} \\
w_{k+1}&w_k\\
\end{array}\right ] \equiv \left[\begin{array}{cr}
qb+a &b \\
b & a\\
\end{array}\right ].
 \eqno(2.2)$$
In the case gcd$\,(a^2+qab-b^2, m)=1$, all  the matrices in (2.2) are
  invertible modulo $m$. Hence
  $k$ is the smallest positive integer such that $\sigma^k\equiv 1$ modulo $m$.
   \qed

\noindent{\bf Discussion 2.2.}
The smallest $k >0$ such that $\sigma^k\equiv 1$ (mod $m$) is called the order of $\sigma$ modulo $m$.
Lemma 2.1 implies that if gcd$\,(a^2+qab -b^2, m)=1$, then the length of $w(a,b,q,m)$ is the order
 of $\sigma$ modulo $m$ which is also the length of $w(0,1,q,m)$.
 In general, the length of a period $w(a,b,q, m)$ is a divisor
 of the order of $\sigma$ modulo $m$.

Suppose that gcd$\,(a^2+qab-b^2,m)=1$.
Applying Lemma 2.1, the length of $w(a,b,q, m)$ is the order of $\sigma$ modulo $m$.
Denoted by $k(m)$ the order of  $\sigma$ modulo $m$.

\noindent {\bf Lemma 2.3.} {\em
 Suppose  gcd$\,(m,n)=1$. Then
  $k(mn) = k(m)k(n)/gcd\,(k(m),k(n))$.
  Further,
   \begin{enumerate}
  \item[(i)] Suppose that gcd$\,(q,5)=1$,
  $q\not\equiv 7, 18\,\,(mod\,\,25)$. Then
  $k(5^{n+1}) = 5k(5^n)$ for all $n \ge 1$.
 \item[(ii)]  If $q\equiv 1,4$ $(mod\,\, 5)$, then $k(5) = 20$.
   If $q\equiv 2,3$  $(mod\,\, 5)$, then $k(5) = 12$.
  \item[(iii)] Suppose that gcd$\,(q,3)=1, q \not\equiv 4, 5\,\,(mod\,\,9)$. Then
  $k(3^n) = 8\cdot 3^{n-1}$.
  \item[(iv)] If $q\equiv 4, 5\,(mod \,9)$, then $k(3) = k(9) =8$. If
 $q\equiv 7,  18\,(mod \,5^2)$, then $k(5)= k(5^2) =12$.
  \end{enumerate}}

  \noindent {\em Proof.} Recall that $\sigma$ is the matrix in (2.1).
   We shall prove (iii) for the case $q \equiv 1$ (mod 9).
    Direct calculation shows that $\sigma^4\not \equiv 1$ (mod 3) and
     $$\sigma^{8} = I + 3A,\eqno(2.3)$$ where
  $A = (a_{ij}) $ is a two by two matrix such that
   gcd$\,(3, a_{11}, a_{12},a_{21}, a_{22})=1$.
   This implies that the order
    of $\sigma$ modulo 3 is 8.
     One may now apply (2.3) and induction to prove that
      $k(3^n) = 8k(3^{n-1})$.
The remaining cases can be proved similarly.
 \qed

%\noindent {\bf Discussion 2.5.}
%Let $f(r)$ be  the length of a period
% of the Fibonacci numbers modulo $r$. Then
%  \begin{enumerate}
  %  \item[(a)]  $f(2) = 3$, $f(5) = 20$,
 %   \item[(b)] $f(r)|(r-1)$ if $r\equiv \pm 1$ (mod 5),  $f(r)|2(r+1)$ if $r\equiv \pm 2$ (mod 5).
 %  \end{enumerate}
%Burr [3] defined a system he called a  complete Fibonacci system modulo $m$.
%Many of his interesting results are consequences of the study of the system.

\noindent {\bf Definition 2.4.}
 Two   periods  $w(a,b,q, m)$ and $w(c,d, q,m)$ are called
 equivalent if one can be obtained from  the other by a cyclic permutation. The set of
  all inequivalent periods is called a   fundamental  system modulo
  $m$. Denote this system by $FS(m)$.

\noindent {\bf Lemma 2.5.} {\em Suppose that gcd$\,(q,m)=1$.
The total number of terms in $FS(m)$ is $m^2-1$.}

 \noindent {\em Proof.} Since
 gcd$\,(q,m)=1$, one sees easily that
  every period has length at least 3.
  Let $ C = (c_1,c_2,\cdots c_k)$ be a period  of length $k$.
    $C$ gives $k$ adjacent pairs
  $(c_1,c_2), (c_2,c_3), $
  $\cdots, (c_{k-1},c_k)$ and $ (c_k,c_1)$.
   Since $k\ge 3$, the above pairs are all distinct from one another.

  Let $(a,b) $ be a nonzero pair. By our definition of periods,
  $(a,b)$ must appear as adjacent terms in some  periods.
   Since the fundamental
   system consists
   of inequivalent periods, the pair  $(a, b)$
   appears  as adjacent terms exactly once
    in $FS(m)$. Our assertion now follows from  the fact that there are exactly
     $m^2 -1$ nonzero pairs in $\Bbb Z_m \times \Bbb Z_m$
     and that each period  $C$ of length $k$ contributes exactly $k$
      adjacent pairs.\qed

 \noindent {\bf Lemma 2.6.} {\em Let $\{w_r\}$ be given as in $(1.1).$  Then
   $w_{n+2}w_{n}-w_{n+1}^2
   = w_{n}^2 +qw_nw_{n+1}-w_{n+1}^2 = (-1)^n
   (a^2+qab-b^2)
   $  for all $n$.
   In particular, $\pm
   (a^2+qab-b^2)$ is an invariant of $\{w_r(a,b,q)\}$.
   }

\noindent {\em Proof.} The lemma follows  by taking the determinants
 of the matrices given in (2.1).\qed

Let $w_n= w_n(a,b,q)$.
 Suppose that $\{w_n\}$ is residue complete modulo $m$. Then $0\in \{w_n\}$ modulo $m$.
 It follows that
  $w(a,b, q,m) = d\cdot w(0,1, q,m)$ for some $d$. Applying Lemma 2.6,
  $a^2 +qab -b^2 \equiv \pm d^2$ modulo $m$.
   Since $ 1\in d\cdot w(0,1, q,m)$,
   we conclude that gcd$\,(d,m)=1$.
    It follows that gcd$\,(a^2 +qab -b^2, m)=1$. In summary, the following is true.

    \noindent {\bf Lemma 2.7.} {\em The recurrence   $\{ w_n( a,b,q)\} $ modulo $m$ is residue
     complete if and only of $ w(a,b,q,m) = d\cdot w(0,1,q,m)$,
      gcd$\,(a^2 +qab -b^2, m)=1$ and $\{w_n(0,1,q)\}$ modulo $m$ is residue complete.
      Suppose that $\{w_n(a,b,q)\} $  is  residue complete modulo $m$. Then
     $\{w_n(a,b,q)\} $ is residue complete modulo $r$ for any $r|m$.}

The following lemmas give some basic results about $\{w_r(a,b, q)\}$ when $p$ is a
 divisor $q^2+4$.

 \noindent {\bf Lemma 2.8.} {\em Let $p$ be a prime divisor of $q^2+4$.
  Then $\{w(a,b, q)\} $ modulo $p$ is residue complete if and only of
  gcd$\,( a^2+qab-b^2, p)=1$.}

\noindent {\em Proof.} See  Bumby [ 2],  Somer [8], and Webb and Long [9].\qed

 \noindent {\bf Lemma 2.9.} {\em
 Let $p$ be an odd prime.
 Suppose that $p|(q^2+4)$ and gcd$\,(a^2+qab -b^2,p)=1$.
  Then the length of $w(a,b,q,p)$ is $4p$. }

\noindent {\em Proof.} Apply results of  Wyler [10] (see Lemma 2.2 and 2.4 of [6] also). \qed

\section {$p\in \Omega$}
The main purpose of this section is to study $\{w_r(a,b,q)\}$ modulo $p$, where $p$ is a
 divisor of $q^2+ 4$. Note that the  proof of Lemma 3.1 can be applied
  to prove Lemma 4.3.

 \noindent {\bf Lemma 3.1.}
 {\em let $p$ be an odd prime. Suppose that   $p|gcd\,(q^2+4,m)$.
   Suppose further that $k(pm) = pk(m)$ and $\{w_r\}$ modulo $m$ is
  residue complete. Then $\{w_r\}$ modulo $pm$ is residue complete.}

 \noindent {\em Proof.} Set $k(m) = k$ and $a^2+qab -b^2 = D$.
  Since
 $\{w_r\}$ modulo $m$ is
  residue complete, gcd$\,(a^2+qab-b^2, m) =1$ and
 for each $ A \in \Bbb Z_m$,
  $w_n \equiv A$ modulo $m$  for some $n$.
   Since the length of a period of $\{w_r\}$ modulo $m $ is $k$, one has
   $w_n \equiv w_{n+k} \equiv \cdots \equiv w_{n+(p-1)k}\equiv A$ modulo $m$.
    Hence
    $$\{ w_{n}, w_{n+k},\cdots, w_{n+(p-1)k}\} \equiv \{ A + i m\, :\, 0\le i\le p-1\}\,\,\,(\mbox{mod } pm).\eqno(3.1)$$
 Set $w_{n+1} \equiv B$ (mod $m)$. Then
  $$\{ w_{n+1}, w_{n+k+1},\cdots, w_{n+(p-1)k+1}\} \equiv \{ B + j m \,:\, 0\le j\le p-1\}\,\,\,(\mbox{mod } pm).\eqno(3.2)$$
   Our goal is to show that members in (3.1) are distinct from one another modulo $pm$.
Suppose that two members in (3.1) are equal to each other modulo $pm$. Without loss of generality,
  $w_n \equiv  w_{n+(p-1)k}$ modulo $pm$. Then $w_n \equiv w_{n+(p-1)k} \equiv A + im $ modulo $pm$ for some $i\le p-1$.
Since $\pm D = \pm(a^2+qab-b^2)$ is an invariant of $\{w_r(a,b,q)\}$ (Lemma 2.6),
$$
w_{n}^2 -qw_{n+1} w_{n}- w_{n+1}^2
 = \pm D \,,\,\,\,
w_{n+(p-1)k}^2 -qw_{n+(p-1)k+1} w_{n+(p-1)k}- w_{n+(p-1)k+1}^2  = \pm D.\eqno(3.3)$$
Since
 $w_n \equiv w_{n+(p-1)k} \equiv A + im $  modulo $pm$,
 equations in (3.3) take the following alternative forms modulo $pm$
{\small $$
(A+im)^2 + qw_{n+1}(A+ im)-w_{n+1}^2 \equiv \pm D\,,\,\,\,
(A+im)^2 + qw_{n+(p-1)k+1}(A+ im)-w_{n+(p-1)k+1}^2 \equiv \pm D.$$}
It follows that $ Y = w_{n+1}$ and $ w_{n+(p-1)k +1}$  are solutions  of the following equation modulo $pm$.
$$(A+im)^2  +qY(A+ im)-Y^2
  \equiv \pm D.\eqno(3.4)$$
Since $w_{n+1}$ and $ w_{n+(p-1)k +1}$ are members in (3.2),
they take the form $B+mj$. Hence the $j$'s associated with  $w_{n+1}$ and $ w_{n+(p-1)k +1}$
 are solutions for  $y$ of the following equation modulo $pm$.
$$ (A+im)^2
 +q(B+ym)(A+ im)- (B+ym)^2
  \equiv \pm D.\eqno (3.5)$$
Note that $A^2 +qAB-B^2\equiv \pm D$ (mod $m$). This implies that $ A^2 +qAB-B^2  = m T \pm D.$
An easy calculation shows that the left hand side of (3.5) takes the following form.
$$
L
=m^2(-y^2 + qiy +i^2) - my(2B-qA) +mi(qB+2A) + mT \pm D
.\eqno(3.6)$$
Equation (3.6) implies that (3.5) holds if and only if
 $L$ is congruent to $\pm D$ modulo $pm$. Since $p|m$ and  gcd$\,(D,m)=1$,  it is equivalent to
$$ -y(2B-qA) + i(qB+2A) + T \equiv 0 \,\,\,(\mbox{mod } p).\eqno(3.7)$$
However, $2B -qA\not \equiv 0$ (mod $p$) since otherwise
$\pm D \equiv A^2 +qAB-B^2\equiv 4^{-1}(q^2+4)A^2\equiv 0\,\,\,(\mbox{mod } p).$
 A contradiction ($p$ is a divisor of $q^2+4$).
Similarly,  $qB+2A\not \equiv 0$ (mod $p$).
As a consequence,
for each $i$, there exists exactly one $y= j$ such that (3.7) (as well as (3.5)) is true.
Hence for each $A+im$, there is a unique $Y$ of the form $B+ym$ such that (3.4) is true.
  Hence $w_{n+1} \equiv  w_{n+(p-1)k+1}$ modulo $pm$.
 This implies that $(w_{n}, w_{n+1})\equiv  (w_{n+(p-1)k}, w_{n+(p-1)k +1})$ modulo $pm$.
  In particular, the length of a  period of $\{w_r\}$ modulo $pm$  is at most $(p-1)k$. This is
   a contradiction.
   Hence the $p$ members in (3.1) are all distinct from one another modulo $pm$.
Since
$ k(pm) = pk(m)$,
$\{w_r\}$ modulo $m$ is residue complete and every member in $\{w_r\}$ modulo $m$
 has $p$  pre-images in $\{w_r\}$ modulo $pm$, we conclude that
  $\{w_r\}$ modulo $pm$ is residue complete.\qed

%\noindent {\bf Definition 3.2.} Let $k\in \Bbb N$ be fixed. Suppose that $p|(q^2+4)$.
%Let $w(a,b,q,p) = (w_0, w_1, \cdots, $
%$w_{4p-1})$ be a period of length $4p$.
%The period  $w(a,b,q,p) $  is called $k$-complete if for any $i$, the terms
% $w_i, w_{(i+k)_p}, w_{(i+2k)_p} \cdots ,  w_{(i+(p-1)k)_p}$
% are not congruent to one another modulo $p$,
%  where $(i+xk)_p$ is the remainder of $i+xk$ divided by $4p$.

%\noindent {\bf Example 3.3.} A period $w(0,1, 1,5)$ of the Fibonacci numbers modulo 5 takes
% the following form
 % $$
  %   (0,1,1,2,3,0,3,3,1,4,0,4,4,3,2,0,2,2,4,1).\eqno(3.8)$$
   %  One sees  that the period is $4x$-complete as long as $x$ is not a multiple of $5$.
   %   See Example 5.4 for another $4x$-complete period.

\noindent {\bf Lemma 3.2.}
 {\em Let $p$ be an odd prime.
  Suppose that $p|(q^2+4)$ and gcd$\,(p, m)= 1$.
 Suppose further that $k(pm) = pk(m)$ and  $\{w_r= w_r(a,b,q)\}$ modulo $m$ is
  residue complete.
   Then $\{w_r = w_r(a,b,q)\}$ modulo $pm$ is residue complete.}

 \noindent {\em Proof.}
 By Lemma 2.7, we may assume $\{w_r\}=\{w_r(0,1,q)\}$.
  Set $k(m)= k$.
 Since $ pk(m) = k(pm)
  = k(p)k/gcd\,(k(p), k)$ and $k(p) =  4p$ (Lemma 2.3 and 2.9), we conclude that
   gcd$\,(k,4p)=4$.
  Since   $\{w_n\} $ modulo $m$    is residue complete, for each $ A\in \Bbb Z_{m}$,
   there exists some $  w_n $
    such that $A\equiv w_ n$ modulo $m$.
  We now consider the set
     $ X =\{ w_n, w_{n+k}, w_{n+2k},\cdots, w_{n+(p-1)k}\}$. Since the length of
      a period of $\{w_n\}$ modulo $m$ is
       $ k $,
     $w_{n+rk} \equiv A$ (mod $m$) for $r=0,1,2, \cdots p-1.$
     We now consider $X$ modulo $p$.
      Since
     \begin{center}
       gcd$\,(k,4p)=4$, $k(p)=4p$, and $w_{y+4p}\equiv w_y$ (mod $p$)
    \end{center}
    for all $y$,   the set $X$ modulo $p$ takes the  form
      $ X = \{w_{x}, w_{x + 4}, w_{x+8}, w_{x+12}, \cdots , w_{x+4(p-1)}\}$
      for some $x$ in the range $0\le x\le 3$.
    By our results in Appendix $A$,
     $X$ modulo $p$ is $\{0,1,2,\cdots , p-1\}=\Bbb Z_p$.
   As a consequence,  members in $X $ are not congruent to one another modulo $pm$.
        In particular, every member in $\{w_r\}$  modulo $m$ has $p$ pre-images
         in $\{w_r\}$ modulo $pm$. Since $k(pm)=pk(m)$, $\{w_n\}$ modulo $pm$  is residue
          complete.\qed

 \section { $p\in \Delta$}
 $\Delta =\{2,3,5,7\}$. We shall study $\{w_r\}$ modulo $p^e$ $(p\in \Delta)$ in the following
 subsections. The main results can be found in lemmas 4.1-4.4.

\smallskip
\noindent {\bf 4.1.} We study $\{w_r\}$ modulo $2^e$. Note that in the case $q$ is even $2$ is also
 a divisor of $q^2+4$.

 \noindent {\bf Lemma 4.1.} {\em Let $w_n= w_n(0,1,q)$ be given as in $(1.1)$.  Then the following
  holds.
  \begin{enumerate}
 \item[(i)] Suppose that $q$ is odd. Then $\{w_n\}$ modulo $4$ is residue complete and
 $\{w_n\}$ modulo $2^n$  is not residue complete for all $n \ge 3$.
 \item[(ii)]Suppose that $q$ is even. Then $\{w_n\}$ modulo $2$ is residue complete and
 $\{w_n\}$ modulo $2^n$  is not residue complete for all $n \ge 2$.
 \end{enumerate}}

\noindent {\em Proof.} The lemma follows easily from (iii) and (iv) of the following,
 which can be verified by direct calculation.
(iii) Suppose that $q$ is odd. Then $\{w_n\}$ modulo 4 is residue complete and
 $\{w_n\}$ modulo 8 is not residue complete.
 (iv) Suppose that $q$ is even. Then $\{w_n\}$ modulo 2 is residue complete and
 $\{w_n\}$ modulo 4 is not residue complete.
 \qed

\noindent {\bf 4.2.} We study $\{w_r\}$ modulo $3^e$. Note that $3$ cannot be a divisor of
 $q^2+4$. The idea of the proof of Lemma 4.2 is mainly taken from Burr [3].

\smallskip
\noindent {\bf Lemma 4.2.} {\em
If $3|q$, then $\{w_n\}$ modulo $3^n$ is not residue complete.
Suppose  gcd$\,(a^2+qab-b^2, 3)
 =\,gcd\,(q, 3) =1$.  In the case $q\equiv 4, 5\,\,(mod\,\,9)$,  $w(a,b, q, 3^n)$ is residue complete
    if and only if $n = 1$.
  In the case $ q\not\equiv 4, 5\,\,(mod \,\,9)$,
  $w(a,b, q, 3^n)$ is residue complete for all $n \ge 1$.
}

\noindent {\em Proof.} The first part of the lemma is clear.
In the case $q\equiv  4,5$ (mod 9), by (iv) of Lemma 2.3, $w(a,b, q, 9)$ is not residue complete.
We now assume
  gcd$\,(q(a^2+qab-b^2), 3)=1$ and $q\not\equiv 4, 5$ (mod 9).
   We shall first prove that  $w(0,1,q,m)$ is residue complete.
By Lemma 2.6, $\pm (a^2+qab-b^2)$ is an invariant for $w(a,b,q,m)$. Since the invariant for
 $w(0,1,q,m)$ is $\pm 1$, our assertion is  proved if the following two facts are verified.
 \begin{enumerate}
 \item[(i)] The fundamental system $FS(m)$ has only one period with invariant $\pm 1$ modulo $m$.
 \item[(ii)] For any $a\in \Bbb Z_m$, there exists some $b$ such that $a^2+qab-b^2\equiv 1$
  or $-1$ modulo $m$.
  \end{enumerate}

 \noindent {\bf Proof of (i) for $m=3^n$}.
 Applying Lemma 2.3, the period $ C = w(0,1, q, 3^n)$ has length $8\cdot 3^{n-1}$.
 For each $r_i$ prime to 3 in the range $1\le r_i \le 3^n/2$, define
  $C_i = \{ r_i x\,:\, x \in C\}$. The $C_i$'s are not equivalent to each other as their invariants
  $\pm r_i^2$ are  distinct from each other modulo $m$.
   For each $D\in FS(3^{n-1})$, define
    $3\cdot D= \{3d\,:\, d \in D\}$.
    They clearly form a subset of inequivalent periods for $FS(3^n)$.
     Let $\phi (x)$ be the Euler  function.
     The total number of  terms in
      $\{ C_i \, :\, 1\le i \le 3^n/2  \} \cup
      \{ 3\cdot D\,:\, D\in FS(3^{n-1})\}$
      is
      $8\cdot 3^{n-1} \phi(3^n)/2 + 3^{2(n-1)} = 3^{2n}-1.$
 Hence
 $$   FS(3^{n}) = \{ C_i \, :\, 1\le i \le 3^n/2  \} \cup
      \{ 3\cdot D\,:\, D\in FS(3^{n-1})\}.\eqno(4.1)$$
   Since $x^2 \equiv \pm 1$ (mod $3^n$) if and only if $ x\equiv \pm 1$ (mod $3^n$), we conclude
 that the only period in (4.1) with invariant $\pm 1$  modulo $3^n$ is $ C_1 = C  = w(0,1, q,3^n)$.

 \noindent {\bf Proof of (ii) for $m=3^n$}.
 We consider two cases
(i) $q$ is a multiple of $3$, (ii) gcd$\,(q, 3) =1$. In the first
 case, one sees easily that
 $a^2+qax-x^2\equiv - 1 \,\,(\mbox{mod }3^n)$ is solvable.
  In the second case,
 $a^2+qax-x^2\equiv 1 \,\,(\mbox{mod }3^n)$ is solvable.

   For any $a \in \Bbb Z_m$, by (ii), the invariant of  $w(a,b,q,m)$ is $\pm 1$ for some $b$.
    By (i), $w(a,b,q,m) = w(0,1,q,m)$. In particular. $a \in w(0,1,q,m)$. Hence $w(0,1,q,m)$ is
     residue complete. Consequently, all the $C_i$'s are residue complete. This
      completes the proof of the lemma.
   \qed

\smallskip
\noindent  {\bf 4.3.} We study $\{w_r\}$ modulo $5^e$. Note that $5$ can be a divisor of $q^2 +4$.

\noindent {\bf Lemma 4.3.}
 {\em  Suppose that $5$ is not a divisor of $q(q^2+4)$,  $5|m$, and
  $a^2+qab-b^2 = D\equiv \pm 1\,\,(mod\,\,5)$.
   Suppose further that $k(5m) = 5k(m)$ and $\{w_r\}$ modulo $m$ is
  residue complete. Then $\{w_r\}$ modulo $5m$ is residue complete.}

\noindent {\em Proof.} Our lemma can be proved by applying the proof of Lemma 3.1 step
by step. The only difference is how the fact
 $2B-qA \not\equiv 0$  (mod $p$) is being verified (see (3.7)).
  In Lemma 3.1,  one has $2B-qA \not\equiv  0$ (mod $p$) since otherwise
$\pm D \equiv A^2 +qAB-B^2\equiv 4^{-1}(q^2+4)A^2\equiv 0\,\,\,(\mbox{mod } p)$, while
 in the present proof, one has
$2B-qA \not\equiv 0$ (mod $5$) since otherwise
$\pm 1\equiv \pm D \equiv A^2 +qAB-B^2\equiv 4^{-1}(q^2+4)A^2\,\,\,(\mbox{mod } 5)$,
 where gcd$\,(q,5)=1$.  A contradiction.\qed

\smallskip

\noindent {\bf 4.4}
 Throughout the section $w_n= w_n(0,1,q)$.
 The purpose of this section is to study  $\{w_r\}$
  modulo $7^e$.  Note that 7 cannot be a divisor of $q^2+4$.
The following is clear.
\begin{enumerate}
 \item[(i)]$q^2 +4$ is  not  a square modulo $7$ if and only if  $\{w_r\}$ modulo 7 is  residue complete,
 \item[(ii)] $\{w_r\}$ modulo 49 is not residue complete.

 \end{enumerate}
 Note that in the first case $q\equiv 1,3,4,6$ (mod 7). In summary, the following holds.

 \noindent {\bf Lemma 4.4.} {\em
 The recurrence $\{w_r= w_r(0,1,q)\}$ modulo $7^n$ is residue complete if and only if
  $n=1$ and $q\equiv 1,3,4,6\,\,(mod\,\,7)$.}

\section { The Main Results}

 \noindent {\bf Proposition  5.1.} {\em  Pell numbers modulo $m$ is residue complete
   if and only if $m \in \{2, 3^a, 5^b\}$.}

   \noindent {\em Proof.}
      Recall first that $P_n = w_n(0,1,2)$. Suppose that $\{P_n\}$
    modulo $m$ is residue complete. By  (1.2),  the possible prime divisors of $m$
    are 2, 3, 5,  or 7. Direct calculation shows that  $w(0,1,2,7)$ is not residue complete.
     Hence  $m$ is of the form $2^a3^b5^c$ (Lemma 2.7). Direct calculation shows that
      $\{P_n\}$ modulo $m$ is not residue complete if $m= 6$, 10, 15.
     Hence  $m\in\{2^a, 3^b, 5^c\}$ (Lemma 2.7).

  Direct calculation shows that $\{w(0,1, 2)\}$ modulo 5 is residue complete.
   By Lemma 4.3 and 2.3,  $\{w(0,1, 2)\}$ modulo $5^c$ is residue complete.
   The proposition  now can be proved by applying Lemma 4.1 and 4.2.\qed

  Recall that Pell-Lucas numbers is given by
   $w(2,2, 2)$.  Similar to the above, the following can be proved as well.

  \noindent {\bf Proposition  5.2.} {\em  Pell-Lucas numbers modulo $m$ is residue complete
   if and only if $m =3^n$.}

\noindent {\bf Discussion 5.3.}
 Suppose that $\{w_n(a,b,q) \}$ modulo $m$ is residue complete.
  Then  $\{w_n(0,1,q) \}$ modulo $m$ is residue complete (Lemma 2.7).
  By (1.2), the set of prime divisors of $m$ is a finite set and
   $m$  is of the form $ 2^a 3^b 5^c 7^d\prod p_i^{e_i}$,
   where $p_i|(q^2+4)$.
  By Lemma 4.1,  one has $a\le 2$. By Lemma 4.4, one has  $d\le 1$.
   Similar to Proposition 5.1,
   one may apply lemmas  3.1, 3.2, 4.1-4.4 and some easy calculation
  to determine completely the set of such $m$'s (see Example 5.4).
  Note that the prime 5 is different from  the other primes. We apply Lemma 3.1 and  3.2 if
    $5$ is a divisor of $q^2+4$, and apply Lemma 4.3 if 5 is not a divisor of
     $q^2+4$.

    In particular, results of Burr [3]
    for Fibonacci numbers $(q^2+1 = 5)$  and
     results of Avila and Chen [1]
     for Lucas numbers ($q^2+4=8$) can be obtained by our method.
     It is our duty to point out that although our result is  more general than
      Burr [3], the main idea is mainly
      taken from Burr [3].

\noindent {\bf Example 5.4.} Define $\{w_r\}$ by $w_0=0$, $w_1=1$, $w_n= 3w_{n-1}+w_{n-2}$.
 Then $q^2+4 = 13$.
  A period of $\{w_r\}$ modulo 13 has length 52
  and is given by
  $$(0,1,3,10,7,5,9,6,1,9,2,2,8,
  0,8,11,2,4,1,7,9,8,7,3,3,12,$$
\vspace{-.5cm}
$$  0,12, 10,3,6,8,4,7,12,4,11,11,5,
  0,5,2,11,9,12,6,4,5,6,10,10,1).$$
   Suppose that $\{w_r\}$ modulo $m $ is residue complete. By (1.2),
  $m$ takes the form $2^a3^b5^c7^d13^e$.
   It is clear that  $\{w_r = w_r(0,1,3)\}$ modulo 3 is not residue complete.
    Hence $m$ takes the form $2^a5^c7^d13^e$.
  By Lemma 4.1, $ a\le 2$. By Lemma 4.4, $d\le 1$.
  Further,
 \begin{enumerate}
 \item[(i)] $\{w_r\}$ modulo 4
 (Lemma 4.1), $5^c$  (Lemma 4.3), $7$ (Lemma 4.4), $13^e$
  (Lemma 3.1) are residue complete.
\end{enumerate}
Direct calculation shows that
$k(2)=3$, $k(4)=6$,
 $k(5^n) =12\cdot 5^{n-1}$,
 $k(7)=16$,
 $k(13^n)=  4\cdot 13^n$,
\begin{enumerate}
 \item[(ii)]$\{w_r\}$ modulo 10, 28, 35  are  not residue complete,
 \item[(iii)] $\{w_r\}$ modulo 14, 52 are  residue complete.
\end{enumerate}
Applying Lemma 3.1 and 3.2, we conclude that $\{w_r\}$ modulo $m $ is residue complete if and only if
$$
m\in\{13^b, 2\cdot 13^b  ,4\cdot 13^b ,5^a 13^b , 7\cdot 13^b , 14\cdot 13^b \}.$$

\section {Appendix A. Residue completeness of subsequences of $\{w_r(0,1,q)\}$}
Set $F_0=0, F_1 =1, F_n=qF_{n-1} +F_{n-2}$ and  $L_0=2, L_1 =q, L_n=qL_{n-1} +L_{n-2}$.
 One can show easily that
 $$F_{n+4} = L_4 F_n -F_{n-4}.\eqno(A.1)$$
 Define four recurrences  $A_n = F_{4n}$, $B_n = F_{4n+1}$, $C_n= F_{4n+2},$ and $D_n= F_{4n+3}$.
  It is clear all four recurrences satisfy the same recurrence
  $$w_n = L_4 w_{n-1}-w_{n-2}.\eqno(A.2)$$
Let $p$ be an odd prime divisor of
   $q^2+4$. Applying Lemma 2.8 and 2.9, $\{F_n\}$ modulo $p$ is residue complete and a period of $\{F_n\}$
   modulo $p$ has length $4p$. Since $p$ is a divisor of $q^2+4$. $L_4 \equiv 2$ (mod $p$).
    Hence $(A.2)$ takes the following form.
    $$w_n = 2 w_{n-1}-w_{n-2}\,\,\,(mod\,\,p).\eqno(A.3)$$
    Consider $\{A_n\}$ modulo $p$. Apply the recurrence relation $(A.3)$, our sequence modulo $p$ takes the following form.
    $$A_0=0, A_1 \equiv -2q, A_2\equiv -4q, \cdots, A_n \equiv -2nq\,,\cdots  \cdots \eqno(A.4)$$
It is clear that $\{A_n\}
= \{ F_0, F_{4}, F_{8},\cdots , F_{4(p-1)}\}=\Bbb Z_p$
 modulo $p$. Similarly, one can show $\{B_n\} =\{C_n\}= \{D_n\}=\Bbb Z_p$ modulo $p$.

\section {Appendix B. $u_n=qu_{n-1}-u_{n-2}$}

Let $u_0=a, u_1=b, u_n=qu_{n-1}-u_{n-2}$. Similar to (2.1), one has the following matrix form.
  $$
\mbox{Let } \tau = \left[\begin{array}{rr}
q &1 \\
-1 &0\\
\end{array}\right ].
 \mbox{  Then }
  \left[\begin{array}{cr}
qb-a &b \\
b & a\\
\end{array}\right ]
\tau^n
=
\left[\begin{array}{lr}
u_{n+2} &u_{n+1} \\
u_{n+1}&u_n\\
\end{array}\right ].
 \eqno(B.1)$$

We shall first give an  alternative proof of the following fact
 (see lemmas 2.3 and 2.4 of [6]).

\noindent {\bf Lemma B1.} {\em Let $p$ be an odd prime and let $k(p)$ be the length of
 a period of $\{u_n\}$ modulo $p$.
 Suppose that gcd$\,(p, -a^2+qab-b^2)=1$.
 Then} (i) $k(p)|(p-1)$ {\em if $x^2 -qx+1\equiv0$ has
  two distinct roots in $\Bbb Z_p$,} (ii)  $k(p)|(p+1)$ {\em if $x^2 -qx+1$ is
   irreducible over
 $\Bbb Z_p$.}

 \noindent
 {\em Proof.}  Since gcd$\,(p, -a^2+qab-b^2)=1$, $k(p)$ is the order of
  $\tau$ modulo $p$. Note that $\tau$ is similar to a diagonal matrix
   diag$\,( \lambda_1, \lambda_2)$, where the $\lambda_i$'s are the
   eigenvalues of $\tau$. In the case $x^2 -qx+1\equiv0$ has
  two distinct roots in $\Bbb Z_p$, $\lambda_i \in \Bbb Z_p^{\times}$. Since $\Bbb Z_p^{\times}$
   is of order $p-1$,
   it follows that the order of $\tau$ is a divisor of $p-1$.
    In the case  $x^2 -qx+1$ is
   irreducible over $\Bbb Z_p$,  since the norm of the $\lambda_i$'s is 1,
   the $\lambda_i$'s are in the kernel
    of the surjective norm map $\mathcal N \,:\, GF(p^2)^{\times} \to \Bbb Z_p^{\times}$,
     where $GF(p^2)$ is a finite field of $p^2$ elements (see Chap VIII, Ex 28 of [4]).
     Since the
      kernel of $\mathcal N$ is  a group of order $p+1$, the order of $\tau$ is a
       divisor of $p+1$.\qed

It is clear that periods of length  $p-1$ or $p+1$ cannot be residue complete. Hence
 $\{u_n\}$ modulo $p$ is residue complete only if $x^2 -qx+1 =0$ has repeated roots.
  Equivalently,  $  p$ is a divisor of $q^2 -4$.
  Let $P$ be the set of all prime divisors of $q^2 -4$. The above argument implies that
   $\{u_n\}$ modulo $m$ is residue complete only of $m$ is of the
    form
    $$\prod _{p_i \in P} p_i^{e_i}.\eqno(B.2)$$
 The sequence $\{u_n\}$ is said to be uniformly distributed modulo $m$ if every element $a$ of
  $\Bbb Z_m$ occurs in a period of $\{u_n\}$ equally often.
  We now state without proof the main result of Bumby [2] (for the sequence $\{u_n\}$).

  \noindent {\bf Theorem B2.} (Bumby [2]) {\em The sequence
   $\{u_n\}$ is uniformly distributed modulo $m$ if and only if
    it is uniformly distributed modulo all prime powers of $m$.
     The sequence  $\{u_n\}$ modulo $p\,(a\, prime)$ is uniformly distributed if and only if
      \begin{enumerate}
  \item[(i)] gcd$\,(p, b-qa/2) =1$ $(if \, q=\pm 2)$, or
   \item[(ii)]gcd$\,(p, \alpha(b-a \alpha))=1$, gcd$\,(p,\alpha -\alpha')\ne 1$ $(if \, q \ne \pm 2)$,
  \end{enumerate}
  where $\alpha$ and $\alpha'$ are roots of $x^2-qx+1 =0$ in $\Bbb C$.
  Equivalently, $p|(q^2-4)$ and gcd$\,(p, -a^2+qab-b^2)=1$. Further,
  If $\{u_n\}$ modulo $p$ is uniformly distributed, then $\{u_n\}$ is
   uniformly distributed  modulo $p^h$ with $h>1$ if and only if }
   (a) $p > 3$, (b) $p=3$, $and $ $q^2 \not\equiv 1\,\,(mod\,\,9)$, or
   (c) $p=2$, $q\equiv 2 \,\,(mod\,\,4)$.

In case (b)
of the theorem, one can show easily that $\{u_n\}$ modulo 9 is not residue complete
 if $q^2\equiv 1$ (mod 9).
 In case  (c) of the theorem, one can show easily that $\{u_n\}$ modulo $4$
  is not residue complete if $q\equiv \pm 1$ (mod 4).
  As a consequence, we have the following corollary of Bumby's Theorem.

 \noindent {\bf Corollary B2.} (Bumby [2]) {\em The sequence
   $\{u_n\}$ is residue complete  modulo $m$ if and only if
    it is residue complete modulo all prime powers of $m$.
     The sequence  $\{u_n\}$ modulo $p\,(a\, prime)$ is residue complete if and only if
    $p|(q^2-4)$ and gcd$\,(p, -a^2+qab-b^2)=1$. Further,
  If $\{u_n\}$ modulo $p$ is residue complete, then $\{u_n\}$ is
   residue complete  modulo $p^h$ with $h>1$ if and only if }
   (a) $p >3$, (b) $p=3$, $and $ $q^2 \not\equiv 1\,\,(mod\,\,9)$, or
   (c) $p=2$, $q\equiv  2\,\,(mod\,\,4)$.

\bigskip

%\noindent f-6-1.tex
\end{document}